\documentclass[12pt,reqno]{amsart}
\usepackage{amsmath,amsthm,amssymb,amsfonts,cite}
\usepackage[mathscr]{eucal}

\usepackage[colorlinks=true, linkcolor=red, citecolor=blue]{hyperref}
\newtheorem{prp}{Proposition}[section]
\newtheorem{thm}[prp]{Theorem}
\newtheorem{lmm}[prp]{Lemma}
\newtheorem{col}[prp]{Corollary}
\theoremstyle{definition}
\newtheorem{dfn}[prp]{Definition}
\newtheorem{rem}[prp]{Remark}

\large\normalsize
\makeatletter
\renewcommand\@seccntformat[1]{\textbf{\@nameuse{the#1}.\hspace{.5em}}}
\renewcommand\section{\@startsection {section}{1}{\z@}%
                                   {-3.5ex \@plus -1ex \@minus -.2ex}%
                                   {2.3ex \@plus.2ex}%
                                   {\normalfont\bfseries\centering}}
\makeatother

\newcommand{\abs}[1]{\lvert#1\rvert}
\newcommand{\norm}[1]{\lVert#1\rVert}
\newcommand{\smatrix}[1]{\left(\!\begin{smallmatrix}#1\end{smallmatrix}\!\right)}
\DeclareMathOperator{\rank}{rank}

\DeclareMathOperator{\cotan}{cotan}
\DeclareMathOperator{\cotanh}{cotanh}
\DeclareMathOperator{\st}{Sin}
\DeclareMathOperator{\ct}{Cos}
\DeclareMathOperator{\ttg}{Tan}
\DeclareMathOperator{\tco}{Cotan}
\DeclareMathOperator{\sh}{Sinh}
\DeclareMathOperator{\ch}{Cosh}
\DeclareMathOperator{\htg}{Tanh}
\DeclareMathOperator{\hco}{Cotanh}
\DeclareMathOperator{\sgn}{sgn}
\newcommand{\Z}{{\mathbb Z}}
\newcommand{\SZ}{\scriptscriptstyle{\Z}}
\newcommand{\N}{{\mathbb N}}

\newcommand{\di}{[0,N]_{\SZ}}
\newcommand{\dii}{[0,N+1]_{\SZ}}
\newcommand{\kn}{{k_{0}}}
\newcommand{\kj}{{k+1}}
\newcommand{\ma}{\mathcal{A}}
\newcommand{\mb}{\mathcal{B}}
\newcommand{\mc}{\mathcal{C}}
\newcommand{\md}{\mathcal{D}}
\newcommand{\ms}{\mathcal{S}}
\newcommand{\mj}{\mathcal{J}}
\newcommand{\Zm}{Z=\smatrix{X \\U}}
\newcommand{\xu}{\smatrix{X \\U}}
\newcommand{\Zsm}{\widetilde{Z}=\smatrix{\widetilde{X} \\ \widetilde{U}}}
\newcommand{\Zpm}{\widehat{Z}=\smatrix{\widehat{X} \\ \widehat{U}}}
\newcommand{\tp}{\mathcal{P}}
\newcommand{\tpo}{\tp^{(1)}}

\newcommand{\tq}{\mathcal{Q}}
\newcommand{\tqo}{\tq^{(1)}}

\newcommand{\hp}{\mathscr{P}}
\newcommand{\hpo}{\hp^{(1)}}

\newcommand{\hq}{\mathscr{Q}}
\newcommand{\hqo}{\hq^{(1)}}

\newcommand{\tpi}{\tp^{(i)}_k}
\newcommand{\tqi}{\tq^{(i)}_k}
\newcommand{\sti}{\st^{(i)}}
\newcommand{\sto}{\st^{(1)}}
\newcommand{\stt}{\st^{(2)}}
\newcommand{\stp}{\st^{\!+}}
\newcommand{\stm}{\st^{\!-}}
\newcommand{\cto}{\ct^{(1)}}
\newcommand{\ctt}{\ct^{(2)}}
\newcommand{\cti}{\ct^{(i)}}
\newcommand{\ctp}{\ct^{\!+}}
\newcommand{\ctm}{\ct^{\!-}}
\newcommand{\ttgo}{\ttg^{(1)}}
\newcommand{\ttgt}{\ttg^{(2)}}
\newcommand{\ttgp}{\ttg^{\!+}}
\newcommand{\ttgm}{\ttg^{\!-}}
\newcommand{\tcoo}{\tco^{(1)}}
\newcommand{\tcot}{\tco^{(2)}}
\newcommand{\tcop}{\tco^{\!+}}
\newcommand{\tcom}{\tco^{\!-}}
\newcommand{\sho}{\sh^{(1)}}
\newcommand{\sht}{\sh^{(2)}}
\newcommand{\shp}{\sh^{\!+}}
\newcommand{\shm}{\sh^{\!-}}
\newcommand{\cho}{\ch^{(1)}}
\newcommand{\cht}{\ch^{(2)}}
\newcommand{\chp}{\ch^{\!+}}
\newcommand{\chm}{\ch^{\!-}}
\newcommand{\htgo}{\htg^{(1)}}
\newcommand{\htgt}{\htg^{(2)}}
\newcommand{\htgp}{\htg^{\!+}}
\newcommand{\htgm}{\htg^{\!-}}
\newcommand{\hcoo}{\hco^{(1)}}
\newcommand{\hcot}{\hco^{(2)}}
\newcommand{\hcop}{\hco^{\!+}}
\newcommand{\hcom}{\hco^{\!-}}

\begin{document}


\newcommand{\mytitle}{\uppercase{\bf
Discrete trigonometric\\ and hyperbolic systems:\\[4mm] An overview
}}
\newcommand{\runninghead}{
Discrete trigonometric and hyperbolic systems
}

\newcommand{\myauthor}{Petr Zem\'{a}nek}
\newcommand{\myaddress}{Department of Mathematics and Statistics,\\[0.5mm] Faculty of Science, Masaryk University,\\[0.5mm]
Kotl\'{a}\v{r}sk\'{a} 2, CZ-61137 Brno, Czech Republic}
\newcommand{\myemail}{E-mail: zemanekp@math.muni.cz}

\newcommand{\mykeywords}{{\it Key words and phrases.}
  Discrete symplectic system; Trigonometric system; Hyperbolic system; Trigonometric function; Hyperbolic function.}
\newcommand{\mysubjclass}{2000 {\it Mathematics Subject Classification.}
  39A12; 39-06; 26D05; 33B10.}
\renewcommand{\thefootnote}{}

\thispagestyle{empty}

\begin{center}
  \vspace*{10mm}
  {\Huge\mytitle} \\[15mm]
  {\large\sc\myauthor} \\[3mm]
  {\rm\myaddress \\[1mm] \myemail} \\[10mm]

    \begin{minipage}{15cm}
    {\small
    {\sc\bf Abstract.}
    In this paper we present an overview of results for discrete trigonometric and hyperbolic systems. These systems are discrete analogues of trigonometric and hyperbolic linear Hamiltonian systems. We show results which can be viewed as discrete $n$-dimensional extensions of scalar continuous trigonometric and hyperbolic formulaes.
   }
  \end{minipage}
  \vfill

  \begin{minipage}[b]{0.36\textwidth}
    \raggedleft {\bf Date (revised final version):}  
  \end{minipage}\hspace{0.01\textwidth}
  \begin{minipage}[t]{0.55\textwidth}
  \today
  \end{minipage}\\[4mm]
  \begin{minipage}[t]{0.36\textwidth}
    \raggedleft {\bf How to cite:}  
  \end{minipage}\hspace{0.01\textwidth}
  \begin{minipage}[t]{0.55\textwidth}
    Ulmer Seminare \"{u}ber Funktionalanalysis und Differentialgleichungen {\bf 14}, 
      pp. 345--359, University of Ulm, Ulm, 2009.
  \end{minipage}\\[4mm]
    \begin{minipage}[t]{0.36\textwidth}
    \raggedleft {\bf License:}  
  \end{minipage}\hspace{0.01\textwidth}
  \begin{minipage}[t]{0.55\textwidth}
    \copyright{} \the\year. This manuscript version is made available under the 
     \href{http://creativecommons.org/licenses/by-nc-nd/4.0/}{CC-BY-NC-ND 4.0 license}.
  \end{minipage}\\[4mm]
  \vfill
\end{center}

\footnotetext{\mysubjclass}
\footnotetext{\mykeywords}


\markboth{{\sc\myauthor}}{{\sc\runninghead}}

\setcounter{page}{0}

\newpage
\renewcommand{\thefootnote}{{\bf\,\alph{footnote}\alph{footnote}\alph{footnote}\,}}

\section{Introduction} \label{sec:1}

In this paper we study the discrete trigonometric and hyperbolic systems and certain properties of their solutions. A motivative factor for this research were results, especially from \cite{oD87:O}, which are known for continuous trigonometric and hyperbolic systems --- but they do not have discrete analogues yet. We present a compact theory for discrete trigonometric and hyperbolic systems. Hence we ``recall'' some results which are known in the literature. In particular, we quote the results for trigonometric systems from \cite{drA97}, where special coefficients were considered (but results remain true for more general matrices as well), and the results for hyperbolic systems from \cite{oD.zP01}. All statements are presented without proofs, because the discrete trigonometric and hyperbolic systems are special cases of their time scales extensions, which were studied in \cite{rH.pZ..:T}, or they can be found in some recent papers.

Continuous trigonometric and hyperbolic systems are special cases of a linear Hamiltonian system, which was paid interest in \cite{wK95:Q,oD91:T,oD92:P,oD96:L,oD98:T,wtR71}. On the other hand, discrete trigonometric and hyperbolic systems are special cases of a discrete symplectic system, which was established and investigated in \cite{drA99,mB.oD97,mB.oD.wK..,oD04:D,oD.wK07:A,oD.wK07:O,rH99:D,rH.vR06:I,yS02,cdA.acP96,pZ07}.

The system of the form
\begin{equation}
\label{CTS}
\tag{CTS}
X'=\tq(t)U \quad \text{and} \quad U'=-\tq(t)X,
\end{equation}
where $t\in[a,b]$, $X(t),\ U(t),$ and $\tq(t)$ are $n\times n$ matrices and additionally the matrix $\tq(t)$ is symmetric for all $t\in[a,b]$, is called a \emph{continuous trigonometric system} (CTS). The origin of the study of scalar and matrix trigonometric functions can be found in \cite{pB06}. Other results were published in \cite{jhB57,oD85:O,oD85:P,oD87:R,gjE66:A,gjE66:O,wtR58}. Discrete scalar and matrix trigonometric functions were studied in \cite{drA97,mB.oD00,mB.oD01:Tr,mB.oD01:Th,pZ07} and more recently in \cite{zD.dS03,oD.sP06}.

Properties of \emph{continuous hyperbolic systems}, i.e., the systems in the form
\begin{equation}
\label{CHS}
\tag{CHS}
X'=\hq(t)U \quad \text{and} \quad U'=\hq(t)X,
\end{equation}
where $t\in[a,b]$, $X(t),\ U(t),$ and $\hq(t)$ are $n\times n$ matrices and additionally the matrix $\hq(t)$ is symmetric for all $t\in[a,b]$, can be found in the work \cite{kF87}. Discrete hyperbolic systems were studied in \cite{oD.zP01,pZ07}.

This paper is organized as follows. In the next section we recall the definition and basic properties of the discrete symplectic systems. In Section~\ref{sec:3} we show all (known and new) results for discrete trigonometric systems. Similar results for discrete hyperbolic systems are presented in Section~\ref{sec:4}. In Section~\ref{sec:5} we show some identities for scalar continuous time trigonometric and hyperbolic functions which we ``cannot'' generalize into the $n$-dimensional discrete case and we explain the reason why it is not possible.

Note, that we consider only real case situation (i.e., for the real-valued coefficients and solutions), but the results hold in the complex domain as well (we only need to replace the transpose of a~matrix by the conjugate transpose, the term ``symmetric'' by ``Hermitian'', and ``orthogonal'' by ``unitary'').
\section{Preliminaries on symplectic system}
\label{sec:2}

The \emph{discrete symplectic system} is the first order linear difference system of the form
\begin{equation}
\label{ss}
z_{k+1}=\ms_k z_k,
\end{equation}
where $k\in \{0,1,2,...,N\}\!:\,=\di$, $N\in \N$, $z_k=\smatrix{x_k\\u_k}\in\mathbb{R}^{2n}$, $x_k,\ u_k\in\mathbb{R}^{n}$, $\ms_k\in\mathbb{R}^{2n\times2n}$, and additionally the matrix $\ms_k$ is symplectic for all $k\in\di$, i.e., $\ms_k^T\mj\ms=\mj$ for all $k\in\di$, where $\mj\!:\,=\smatrix{\phantom{-}0 & I\\-I & 0}$.

Since any symplectic matrix is invertible, it is obvious, that symplectic system \eqref{ss} supplemented with the initial condition $z_{\kn}=z^{(0)}$, where $\kn\in\{0,1,...,N+1\}\!:\,=\dii$ and $z^{(0)}\in\mathbb{R}^{2n}$, has a unique solution. Hence we can consider symplectic system \eqref{ss} as a matrix system
\[
\begin{pmatrix}
X_\kj\\
U_\kj
\end{pmatrix}=
\ms_k
\begin{pmatrix}
X_k\\
U_k
\end{pmatrix},
\]
where $X_k,\ U_k\in\mathbb{R}^{n\times n}$. Moreover, when we use the block notation $\ms_k=\smatrix{\ma_k&\mb_k\\\mc_k&\md_k}$, where $\ma_k,\ \mb_k,\ \mc_k,\ \md_k\in\mathbb{R}^{n\times n}$, we can rewrite symplectic system \eqref{ss} in the form
\begin{equation}
\tag{$\ms$}
\label{s}
X_\kj=\ma_k X_k+\mb_k U_k, \quad \text{and} \quad U_\kj=\mc_k X_k+\md_k U_k.
\end{equation}

The condition for the symplecticity of the matrix $\ms_k$ takes the following (using the block notation) equivalent form
\begin{align}
\ma_k^T\md_k-\mc_k^T\mb_k=I=\md_k^T\ma_k-\mb_k^T\mc_k, \quad \ma_k\md_k^T-\mb_k\mc_k^T=I=\md_k\ma_k^T-\mc_k\mb_k^T,\label{3}\\
\ma_k^T\mc_k-\mc_k^T\ma_k=0=\mb_k^T\md_k-\md_k^T\mc_k, \quad \ma_k\mb_k^T-\mb_k\ma_k^T=0=\mc_k\md_k^T-\md_k\mc_k^T.\label{4}
\end{align}

If $\Zm$ and $\Zsm$ are any matrix solutions of \eqref{s}, then their \emph{Wronskian matrix} is defined on $\dii$ as
\[
W(Z,\widetilde{Z})\!:\,=Z^T\mj \widetilde{Z}=X^T\widetilde{U}-U^T\widetilde{X}.
\]
It is known, see \cite[Remark 1\,(ii)]{mB.oD97}, that the Wronskian matrix is constant on $\dii$. A solution $Z=\smatrix{X\\U}$ of \eqref{s} is said to be a \emph{conjoined solution} of \eqref{s} if $W(Z,Z)\equiv0$, i.e., the matrix $X^TU$ is symmetric at one (and hence for all) $k \in \dii$. Two solutions $Z$ and $\widetilde{Z}$ are called \emph{normalized} if $W(Z,\widetilde{Z})\equiv I$. A solution $Z$ is said to be \emph{basis} if $\rank Z=n$ for all $k\in\dii$. It is known that for any conjoined basis $Z$ there always exists another basis $\widetilde{Z}$ such that $Z$ and $\widetilde{Z}$ are normalized, see \cite[Remark 1\,(ii)]{mB.oD97}.

If the matrix solutions $\Zm$ and $\Zsm$ form normalized conjoined bases the following identities
\begin{align}
X_k^T\widetilde{U}_k-U_k^T\widetilde{X}_k&=I=X_k\widetilde{U}_k^T-\widetilde{X}_kU_k^T,\label{5}\\
X_k\widetilde{X}_k^T-\widetilde{X}_kX_k^T&=0=U_k\widetilde{U}_k^T-\widetilde{U}_kU_k^T\label{6}
\end{align}
are true for all $k\in\dii$. Moreover, for all $k\in\di$ the following identities hold
\begin{align}
X_\kj\widetilde{U}_k^{T}-\widetilde{X}_\kj U_k^{T}&=\ma_k,\label{7}\\
\widetilde{X}_\kj X_k^{T}-X_\kj\widetilde{X}_k^{T}&=\mb_k,\label{8}\\
U_\kj\widetilde{U}_k^{T}-\widetilde{U}_\kj U_k^{T}&=\mc_k,\label{9}\\
\widetilde{U}_\kj X_k^{T}-U_\kj\widetilde{X}_k^{T}&=\md_k. \label{10}
\end{align}
The matrix solution $\Zpm$ of the initial value problem \eqref{s}, $X_{\kn}=0$ and $U_{\kn}=I$, with $\kn\in\dii$, is called the \emph{principal solution at $\kn$}.

\section{Discrete trigonometric system}
\label{sec:3}

In this section we consider the symplectic system \eqref{s} with the matrix
\[
\ms_k=
\begin{pmatrix}
\phantom{-}\tp_k & \tq_k\\
-\tq_k & \tp_k
\end{pmatrix},
\]
where $\tp_k,\ \tq_k\in\mathbb{R}^{n\times n}$ for $k\in\di$. From \eqref{3} and \eqref{4} we get that the matrices $\tp_k$ and $\tq_k$ satisfy
\begin{align}
\tp_k^T\tp_k+\tq_k^T\tq_k=I=\tp_k\tp_k^T+\tq_k\tq_k^T, \label{11}\\
\tp_k^T\tq_k-\tq_k^T\tp_k=0=\tp_k\tq^T_k-\tq_k\tp^T_k. \label{12}
\end{align}

The following definition is in accordance with definition of discrete trigonometric system from \cite{drA97}.
\begin{dfn}[Discrete trigonometric system]
\label{3.1}
The system
\begin{equation}
\tag{DTS}
X_{\kj}=\tp_k X_k+\tq_k U_k, \quad U_{\kj}=-\tq_k X_k+\tp_k U_k, \label{ts}
\end{equation}
where the matrices $\tp_k$ and $\tq_k$ satisfy identities \eqref{11}, \eqref{12} for all $k\in\di$, is called a \emph{discrete trigonometric system}.
\end{dfn}
\begin{rem}
\label{3.2}
By a straightforward calculation we may show that the matrix $\ms_k$ for the trigonometric system satisfies $\mj^T\ms_k\mj=\ms_k$ for all $k\in\di$. Therefore, trigonometric systems are \emph{self\discretionary{-}{-}{-}reciprocal}, see e.g. \cite{oD.rH99}.
\end{rem}

\begin{rem}
\label{3.3}
We justify the terminology ``trigonometric'' system. Let us consider scalar discrete trigonometric system
\begin{equation}
\label{13}
x_\kj=p_kx_k+q_ku_k, \quad \text{and} \quad u_\kj=-q_kx_k+p_ku_k
\end{equation}
with the initial conditions $x_0=0$ and $u_0=1$. The coefficients must satisfy (in order to be trigonometric) $p_k^2+q_k^2=1$ for all $k\in\di$. Hence the motivation is the fact that there exist $\varphi_k\in(0,2\pi]$ such that $\cos \varphi_k=p_k$ and $\sin \varphi_k=q_k$. Thus, the pair
\[
x_k=\sin\left(\sum_{j=0}^{k-1}\varphi_j\right) \quad \text{and} \quad u_k=\cos\left(\sum_{j=0}^{k-1}\varphi_j\right)
\]
is the solution of \eqref{13}.
\end{rem}

By an easy calculation we can verify the following lemma.
\begin{lmm}
\label{3.4}
The pair $\xu$ solves system \eqref{ts} if and only if the pair $\smatrix{\phantom{-}U \\-X}$ solves the same trigonometric system. Equivalently $\smatrix{U \\X}$ solves \eqref{ts} if and only if $\smatrix{-X\\\phantom{-}U}$ solves \eqref{ts}.
\end{lmm}

\begin{dfn}[Discrete matrix-valued trigonometric functions]
\label{3.5}
Let $\kn\in\dii$ be fixed. We define the discrete $n\times n$ matrix-valued functions \emph{sine} (denote $\st_k$) and \emph{cosine} (denote $\ct_k$) by
\[
\st_{k;\kn}\!:\,=X_k, \quad \ct_{k;\kn}\!:\,=U_k,
\]
where the pair $\xu$ is the principal solution at $\kn$ of \eqref{ts}. We suppress the index $\kn$ when $\kn=0$, i.e., we denote $\st_k\!:\,=\st_{k;0}$ and $\ct_k\!:\,=\ct_{k;0}$.
\end{dfn}

\begin{rem}
\label{3.6}
By using Lemma~\ref{3.5}, we may define the matrix-valued trigonometric functions alternatively --- as $\ct_{k;\kn}\!\!:\,=\widetilde{X}_k$ and $\st_{k;\kn}\!\!:\,=-\widetilde{U}_k$, where the pair $\smatrix {\widetilde{X}\\ \widetilde{U}}$ is the solution of system \eqref{ts} with the initial conditions $\widetilde{X}_\kn=I$ and $\widetilde{U}_\kn=0$.
\end{rem}

\begin{rem}
\label{3.6.1}
The solutions of system \eqref{CTS} with $n=1$ take the form $\st(t)=\sin\int_a^t\tq(\tau)\,{\mathrm d}\tau$ and $\ct(t)=\cosh\int_a^t\tq(\tau)\,{\mathrm d}\tau$. Hence, we can see the discrete matrix functions $\st_{k;\kn}$ and $\ct_{k;\kn}$ as $n$-dimensional discrete analogs of the scalar continuous trigonometric functions for $\sin(t-s)$ and $\cos(t-s)$, which are solution of scalar system \eqref{CTS} with $\tq(t)\equiv1$.
\end{rem}

From the initial conditions and from the constancy of the wronskian matrix follows, that $\smatrix{\ct\\ \st}$ and $\smatrix{-\st\\ \phantom{-}\ct}$ form normalized conjoined bases of system \eqref{ts} and the matrix
\[
\Phi_k=
\begin{pmatrix}
\ct_k & -\st_k\\
\st_k & \phantom{-}\ct_k
\end{pmatrix}
\]
is a fundamental matrix of \eqref{ts}. Hence, we can express every solution $\smatrix{X\\U}$ of \eqref{ts} in the form
\[
X_k=\ct_kX_0-\st_kU_0, \quad \text{and} \quad U_k=\st_kX_0+\ct_kU_0
\]
for all $k\in\dii$. This statement corresponds to \cite[Lemma~1]{drA97}.

The following identities are consequences of formulaes \eqref{5}--\eqref{10} and the next two corollaries can be found in \cite[Theorem~5~and~Lemma~7]{drA97}.

\begin{col}
\label{3.7}
For all $k\in\dii$ the identities
\begin{align}
\ct_k^T\ct_k+\st_k^T\st_k&=I=\ct_k\ct_k^T+\st_k\st_k^T,\label{14}\\
\ct_k^T\st_k-\st_k^T\ct_k&=0=\ct_k\st_k^T-\st_k\ct_k^T \label{15}
\end{align}
hold, while for all $k\in\di$ we have the identities
\begin{align}
\ct_\kj\ct_k^T+\st_\kj\st_k^T=\tp_k,\label{16}\\
\ct_\kj\st_k^T-\st_\kj\ct_k^T=\tq_k.\label{17}
\end{align}
\end{col}

From formula \eqref{14} we get the following identity, which is a discrete matrix analogue of the fundamental formula $\cos^2 (t)+\sin^2 (t)=1$ for scalar continuous time trigonometric functions. Here $\norm{\cdot}$ is the usual \emph{Frobenius norm}, i.e., $\norm{V}_F=\big(\sum_{i,j=1}^{n}v_{ij}^{2}\big)^{\frac{1}{2}}$, see \cite[pg.~346]{dsB05}.
\begin{col}
\label{3.8}
For all $k \in \dii$ we have the identity
\begin{equation}
\label{18}
\norm{\ct}_{F}^{2}+\norm{\st}_{F}^{2}=n.
\end{equation}
\end{col}


\begin{rem}
\label{3.10}
If the matrix $\ct$ is invertible for some $k\in\dii$, then we can calculate from \eqref{14} and \eqref{15} that
\begin{equation*}
\ct_k^{-1}=\ct_k^{T}+\st_k^{T}\ct_k^{T-1}\st_k^{T}.
\end{equation*}
Similarly,
\begin{equation*}
\st^{-1}_k=\st_k^{T}+\ct_k^{T}\st_k^{T-1}\ct_k^{T}
\end{equation*}
holds if the matrix $\st_k$ is invertible for some $k\in\dii$.
\end{rem}

The following theorem is a discrete analogue of \cite[Theorem~1.1]{gjE66:O} and is contained in \cite[Theorem~9~and~Corollary~10]{drA97}.
\begin{thm}
\label{3.11}
For $k,l\in\dii$ we have
\begin{align}
\st_{k;l}&=\st_k\ct^{T}_l-\ct_k\st^{T}_l, \label{21}\\
\ct_{k;l}&=\ct_k\ct^{T}_l+\st_k\st^{T}_l,\label{22} \\
\st_k&=\st_{k;l}\ct_l+\ct_{k;l}\st_l,\label{23}\\
\ct_k&=\ct_{k;l}\ct_l-\st_{k;l}\st_l. \label{24}
\end{align}
\end{thm}

\begin{rem}
\label{3.12}
With respect to Remark \ref{3.6.1}, identities \eqref{21} and \eqref{22} are discrete matrix extensions of the formulae
\begin{align*}
\sin(t-s)=\sin(t)\cos(s)-\cos(t)\sin(s),\\
\cos(t-s)=\cos(t)\cos(s)+\sin(t)\sin(s).
\end{align*}
\end{rem}

Interchanging the parameters $k$ and $l$ in \eqref{21} and \eqref{22} yields to the following identities, which are $n$-dimensional discrete analogues of the statement about the parity for the scalar goniometric functions. This statement can be found in \cite[Corollary~11]{drA97}.
\begin{col}
\label{3.13}
Let $k,l\in\dii$. Then
\begin{equation}
\label{26}
\st_{k;l}=-\st^{T}_{l;k} \quad \text{and} \quad \ct_{k;l}=\ct^{T}_{l;k}.
\end{equation}
\end{col}

In Theorem~\ref{3.11} we generalized the sum and difference formulaes for solutions of \emph{one trigonometric system} with \emph{different initial conditions}. Now, we shall deal with \emph{two trigonometric systems} and their solutions with the \emph{same initial conditions}. In the continuous case it was done in \cite[Theorem~1]{oD87:O}.

Consider the following two discrete trigonometric systems
\begin{equation}
\label{25}
X_{k+1}=\tpi\, X_k+\tqi\, U_k \quad \text{and} \quad U_{k+1}=-\tqi\, X_k+\tpi\, U_k,
\end{equation}
with initial conditions $X^{(i)}_0=0$ and $U^{(i)}_0=I$, where $i=1,2$. We denote by $\sti_k$ and $\cti_k$ the corresponding solutions of system~\eqref{25} by Definition~\ref{3.5}. Note, that for simplicity we use the following notation $(\tqi)^T=\tq^{(i)T}_k$, $(\tpi)^T=\tp^{(i)T}_k$, $(\sti_k)^T=\st^{(i)T}_k$, $(\cti_k)^T=\ct^{(i)T}_k$, $(\stp_k)^T=\st^{\!+T}_k$, $(\ctp_k)^T=\ct^{\!+T}_k$, $(\stm_k)^T=\st^{\!-T}_k$, and $(\ctm_k)^T=\ct^{\!-T}_k$. Now, we put
\begin{align}
\stp_k&\!:\,=\sto_k\ct^{(2)T}_k+\cto_k\st^{(2)T}_k,\label{27}\\
\ctp_k&\!:\,=\cto_k\ct^{(2)T}_k-\sto_k\st^{(2)T}_k,\label{28}\\
\stm_k&\!:\,=\sto_k\ct^{(2)T}_k-\cto_k\st^{(2)T}_k,\label{29}\\
\ctm_k&\!:\,=\cto_k\ct^{(2)T}_k+\sto_k\st^{(2)T}_k.\label{30}
\end{align}

\begin{thm}
\label{3.14}
Let the matrices $\tpi$ and $\tqi$, $i=1,2$, satisfy \eqref{11} and \eqref{12}. The pair $\stp_k$ and $\ctp_k$ solves the system
\begin{align*}
X_{k+1}&=\phantom{-}\tpo_k(X_k\tp^{(2)T}_k+U_k\tq^{(2)T}_k)+\tqo_k(-X_k\tq^{(2)T}_k+U_k\tp^{(2)T}_k),\\
U_{k+1}&=         - \tqo_k(X_k\tp^{(2)T}_k+U_k\tq^{(2)T}_k)+\tpo_k(-X_k\tq^{(2)T}_k+U_k\tp^{(2)T}_k),
\end{align*}
with the initial conditions $X_0=0$ and $U_0=I$. The pair $\stm_k$ and $\ctm_k$ solves the system
\begin{align*}
X_{k+1}&=\phantom{-}\tpo_k(X_k\tp^{(2)T}_k-U_k\tq^{(2)T}_k)+\tqo_k(X_k\tq^{(2)T}_k+U_k\tp^{(2)T}_k),\\
U_{k+1}&=         - \tqo_k(X_k\tp^{(2)T}_k-U_k\tq^{(2)T}_k)+\tpo_k(X_k\tq^{(2)T}_k+U_k\tp^{(2)T}_k),
\end{align*}
with the initial conditions $X_0=0$ and $U_0=I$. Moreover, for all $k\in\dii
$ we have
\begin{align}
\stp_k\st^{\!+T}_k+\ctp_k\ct^{\!+T}_k=I&=\st^{\!+T}_k\stp_k+\ct^{\!+T}_k\ctp_k, \label{31}\\
\stm_k\st^{\!-T}_k+\ctm_k\ct^{\!-T}_k=I&=\st^{\!-T}_k\stm_k+\ct^{\!-T}_k\ctm_k, \label{32}\\
\stp_k\ct^{\!+T}_k-\ctp_k\st^{\!+T}_k=0&=\st^{\!+T}_k\ctp_k-\ct^{\!+T}_k\stp_k, \label{33}\\
\stm_k\ct^{\!-T}_k-\ctm_k\st^{\!-T}_k=0&=\st^{\!-T}_k\ctm_k-\ct^{\!-T}_k\stm_k. \label{34}
\end{align}
\end{thm}

\begin{rem}
\label{3.15}
Although the above two systems are visually of the same form as the discrete trigonometric system \eqref{ts} and identities \eqref{31} and \eqref{33} for the pair $\stp_k$ and $\ctp_k$ and identities \eqref{32} and \eqref{34} for $\stm_k$ and $\ctm_k$ appear to be like the properties of normalized conjoined bases of system \eqref{s}, the solutions $\smatrix{\stp\\ \ctp}$ and $\smatrix{\stm\\ \ctm}$ are not conjoined bases of their corresponding systems, because these systems are not symplectic.
\end{rem}

When the two systems in \eqref{25} are the same, we obtain the following corollary. It is discrete analogue of \cite[Theorem~1.1]{gjE66:A}.
\begin{col}
\label{3.16}
Assume that the matrices $\tp_k$ and $\tq_k$ satisfy \eqref{11} and \eqref{12}. Then the system
\begin{align*}
X_{k+1}&=\phantom{-}\tp_k(X_k\tp_k^T+U_k\tq_k^T)+\tq_k(-X_k\tq_k^T+U_k\tp_k^T),\\
U_{k+1}&=         - \tq_k(X_k\tp_k^T+U_k\tq_k^T)+\tp_k(-X_k\tq_k^T+U_k\tp_k^T)
\end{align*}
with initial conditions $X_0=0$ and $U_0=I$ possesses the solution
\[
X_k=2\st_k\ct_k^T, \quad \text{and}\quad U_k=\ct_k\ct_k^T-\st_k\st_k^T,
\]
where the functions $\st_k$ and $\ct_k$ are the matrix functions in Definition~\ref{3.5}. Moreover, the above matrices $X$ and $U$ commute, i.e., $X_kU_k=U_kX_k$.
\end{col}

\begin{rem}
\label{3.17}
The previous corollary can be viewed as the discrete $n-$dimensional analogy of the double angle formulae for scalar continuous time trigonometric functions
\[
\sin(2t)=2\sin(t)\cos(t) \quad \text{and} \quad \cos(2t)=\cos^{2}(t)-\sin^{2}(t).
\]
\end{rem}

\begin{col}
\label{3.18}
For all $k \in \dii$ the following identities hold
\begin{align}
\sto_k\st_k^{(2)T}&=\frac{1}{2}\,(\ctm_k-\ctp_k),\label{35}\\
\cto_k\ct_k^{(2)T}&=\frac{1}{2}\,(\ctm_k+\ctp_k),\label{36}\\
\sto_k\ct_k^{(2)T}&=\frac{1}{2}\,(\stm_k+\stp_k).\label{37}
\end{align}
\end{col}

\begin{rem}
\label{3.19}
Identities \eqref{35}--\eqref{37} are discrete $n$-dimensional generalizations of
\begin{align*}
\sin(t)\sin(s)&=\frac{1}{2}\,[\cos(t-s)-\cos(t+s)],\\
\cos(t)\cos(s)&=\frac{1}{2}\,[\cos(t-s)+\cos(t+s)],\\
\sin(t)\cos(s)&=\frac{1}{2}\,[\sin(t-s)+\sin(t+s)].
\end{align*}
\end{rem}

According to the definition in \cite[pg.~42]{drA97} we introduce the discrete $n$-dimensional analogues of scalar trigonometric functions tangent and cotangent.
\begin{dfn}
\label{3.20}
Whenever the matrix $\ct_k$ is invertible we define the discrete matrix-valued function \emph{tangent} (we write $\ttg$) by
\[
\ttg_k\!:\,=\ct^{-1}_k\,\st_k.
\]
Whenever the matrix $\st_k$ is invertible we define the discrete matrix-valued function \emph{cotangent} (we write $\tco$) by
\[
\tco_k\!:\,=\st^{-1}_k\,\ct_k.
\]
\end{dfn}

In the following two theorems we show properties of discrete matrix-valued functions tangent and cotangent, which were published in \cite[Corollary~6~and~Lemma~12]{drA97}.
\begin{thm}
\label{3.21}
Whenever $\ttg_k$ is defined we have
\begin{align}
\ttg^{T}_k&=\ttg_k,\label{40}\\
\ct^{-1}_k\ct^{T-1}_k&-\ttg^{2}_k=I. \label{41}
\end{align}
Moreover, if the matrices $\ct_k$ and $\ct_\kj$ are invertible, then
\begin{equation}
\label{42}
\Delta\ttg_k=\ct_\kj^{-1}\tq_k\,\ct^{T-1}_k.
\end{equation}
\end{thm}

\begin{rem}
\label{3.21.1}
In the scalar case $n=1$ identity \eqref{40} is trivial. In the scalar continuous time case \eqref{CTS} with $\tq(t)=q(t)$ identity \eqref{41} takes the form
\[
\frac{1}{\cos^{2}(s)}-\tan^{2}(s)=1 \quad \text{with \ $s=\int_a^t\!\!\!q(\tau)\,{\mathrm d}\tau\neq\frac{k\pi}{2}$.}
\]
Finally, for the above $s$ identity \eqref{42} takes the form
\[
\left(\tan\!\int_a^t\!\!\!q(\tau)\,{\mathrm d}\tau\right)'=\frac{q(t)}{\cos^{2}\int_a^t q(\tau)\,{\mathrm d}\tau}.
\]
\end{rem}

\begin{thm}
\label{3.22}
Whenever $\tco_k$ is defined we get
\begin{align}
\tco^{T}_k&=\tco_k,\label{43}\\
\st^{-1}_k\,\st^{T-1}_k&-\tco^{2}_k=I. \label{44}
\end{align}
Moreover, if the matrices $\st_k$ and $\st_\kj$ are invertible, then
\begin{equation}
\label{45.1}
\Delta\tco_k=-\st_\kj^{-1}\tq_k\,\st^{T-1}_k.
\end{equation}
\end{thm}
\begin{rem}
\label{3.22.1}
In the scalar case $n=1$ identity \eqref{43} is trivial. In the scalar continuous time case \eqref{CTS} with $\tq(t)=q(t)$ identity \eqref{44} reads as
\[
\frac{1}{\sin^{2}(s)}-\cotan^{2}(s)=1,
\]
where $s\neq k\pi$ is from Remark \ref{3.21.1}. And for these values of $s$ identity \eqref{45.1} reduces to
\[
\left(\cotan\!\int_a^t\!\!\!q(\tau)\,{\mathrm d}\tau\right)'=\frac{-q(t)}{\sin^{2}\int_a^t q(\tau)\,{\mathrm d}\tau}.
\]
\end{rem}

Next, similarly to the definitions of the discrete matrix functions $\sti_k$, $\cti_k$, $\st^{\pm}_k$, and $\ct^{\pm}_k$ from \eqref{27}--\eqref{30} we define the following functions
\begin{align*}
\ttg^{(i)}_k&\!:\,=(\cti_k)^{-1}\,\sti_k,&\quad\tco^{(i)}_k&\!:\,=(\sti_k)^{-1}\,\cti_k,\\
\ttgp_k&\!:\,=(\ctp_k)^{-1}\stp_k,&\quad\tcop_k&\!:\,=(\stp_k)^{-1}\ctp_k,\\
\ttgm_k&\!:\,=(\ctm_k)^{-1}\stm_k,&\quad\tcom_k&\!:\,=(\stm_k)^{-1}\ctm_k.
\end{align*}

\begin{rem}
\label{3.23}
Of course, the matrix-valued functions $\ttg^{\pm}_k$ and $\tco^{\pm}_k$ have similar properties as the functions $\ttg_k$ and $\tco_k$. In particular, it follows from \eqref{33} and \eqref{34} that $\ttg^{\pm}_k$ and $\tco^{\pm}_k$ are symmetric.
\end{rem}

\begin{thm}
\label{3.24}
For all $k\in\dii$ such that all involved functions exist, the following identities hold
\begin{align}
\ttgo_k+\ttgt_k&=\ttgo_k\,(\tcoo_k+\tcot_k)\,\ttgt_k, \label{49.t}\\
\ttgo_k-\ttgt_k&=\ttgo_k\,(\tcot_k-\tcoo_k)\,\ttgt_k, \label{50.t}\\
\ttgo_k+\ttgt_k&=(\ctt_k)^{-1}\,\st_k^{\!+T}(\cto_k)^{T-1}, \label{45}\\
\ttgo_k-\ttgt_k&=(\ctt_k)^{-1}\,\st_k^{\!-T}(\cto_k)^{T-1}, \label{46}\\
\ttgp_k&=(\ctt_k)^{T-1}\,(I-\ttgo_k\ttgt_k)^{-1}\,(\ttgo_k+\ttgt_k)\,\ct_k^{(2)T}, \label{47}\\
\ttgm_k&=(\ctt_k)^{T-1}\,(I+\ttgo_k\ttgt_k)^{-1}\,(\ttgo_k-\ttgt_k)\,\ct_k^{(2)T}, \label{48}\\
\tcoo_k+\tcot_k&=\tcoo_k\,(\ttgo_k+\ttgt_k)\,\tcot_k, \label{49}\\
\tcoo_k-\tcot_k&=\tcoo_k\,(\ttgt_k-\ttgo_k)\,\tcot_k, \label{50}
\end{align}
\begin{align}
\tcoo_k+\tcot_k&=(\stt_k)^{-1}\,\st_k^{\!+T}(\sto_k)^{T-1}, \label{51}\\
\tcoo_k-\tcot_k&=-(\stt_k)^{-1}\,\st_k^{\!-T}(\sto_k)^{T-1}. \label{52}\\
\tcop_k&=(\stt_k)^{T-1}\,(\tcoo_k+\tcot_k)^{-1}\, \notag \\
  &\ \hspace*{45mm} \times(\tcoo_k\tcot_k-I)\,\st_k^{(2)T}, \label{47.c}\\
\tcom_k&=(\stt_k)^{T-1}\,(\tcot_k-\tcoo_k)^{-1}\, \notag \\
  &\ \hspace*{45mm} \times(\tcoo_k\tcot_k+I)\,\st_k^{(2)T}, \label{48.c}
\end{align}
\end{thm}

\begin{rem}
\label{P3.11}
In the scalar continuous time case \eqref{CTS} with $\tq(t)\equiv 1$ take the identities \eqref{49.t} and \eqref{50.t} the following form
\begin{align*}
\tan(t)+\tan(s)&=[\cotan(t)+\cotan(s)]\tan(t)\tan(s), \\
\tan(t)-\tan(s)&=[\cotan(s)-\cotan(t)]\tan(t)\tan(s),
\end{align*}
respectively, identities \eqref{45} and \eqref{46} correspond to
\[
\tan(t)+\tan(s)=\frac{\sin(t+s)}{\cos(t)\cos(s)}, \quad
\tan(t)-\tan(s)=\frac{\sin(t-s)}{\cos(t)\cos(s)},
\]
respectively, identities \eqref{47} and \eqref{48} have the form
\[
\tan(t+s)=\frac{\tan(t)+\tan(s)}{1-\tan(t)\tan(s)}, \quad
\tan(t-s)=\frac{\tan(t)-\tan(s)}{1+\tan(t)\tan(s)},
\]
respectively. Although the identities in \eqref{49} and \eqref{50} reduce to
\begin{align*}
\cotan(t)+\cotan(s)&=[\tan(t)+\tan(s)]\cotan(t)\cotan(s), \\
\cotan(t)-\cotan(s)&=[\tan(s)-\tan(t)]\cotan(t)\cotan(s),
\end{align*}
it is common to write them as
\[
\tan(t)\tan(s)=\frac{\tan(t)+\tan(s)}{\cotan(t)+\cotan(s)}=\frac{\tan(s)-\tan(t)}{\cotan(t)-\cotan(s)}.
\]
The identities in \eqref{51} and \eqref{52} are equivalent to
\begin{align*}
\cotan(t)+\cotan(s)&=\frac{\sin(t+s)}{\sin(t)\,\sin(s)}, \\
\cotan(t)-\cotan(s)&=-\frac{\sin(t-s)}{\sin(t)\,\sin(s)}=\frac{\sin(s-t)}{\sin(t)\,\sin(s)},
\end{align*}
respectively. Finally, the identities \eqref{47.c} and \eqref{48.c} are in accordance with
\[
\cotan(t+s)=\frac{\cotan(t)\cotan(s)-1}{\cotan(t)+\cotan(s)}, \quad
\cotan(t-s)=\frac{\cotan(t)\cotan(s)+1}{\cotan(s)-\cotan(t)},
\]
respectively.
\end{rem}
\section{Discrete hyperbolic system}
\label{sec:4}

In this section we define discrete hyperbolic matrix functions and show analogous results as for the trigonometric functions in the previous section. Some of these results are known from \cite{oD.zP01} but some results and identities are new.

In this section we consider the system \eqref{s} with the matrix $\ms_k$ in the form
\[
\ms_k=\begin{pmatrix}
\hp_k & \hq_k\\
\hq_k & \hp_k
\end{pmatrix},
\]
where the coefficients matrices $\hp_k,\ \hq_k\in\mathbb{R}^{n\times n}$ satisfy for all $k\in\di$ the identities
\begin{align}
\hp_k^T\hp_k-\hq_k^T\hq_k=I=\hp_k\hp_k^T-\hq_k\hq_k^T, \label{H1}\\
\hp_k^T\hq_k-\hq_k^T\hp_k=0=\hp_k\hq_k^T-\hq_k\hp_k^T, \label{H2}
\end{align}
see also \cite[identity~(10)]{oD.zP01}.

\begin{rem}
\label{4.1}
It was shown in \cite{oD.zP01} that from identities \eqref{H1} and \eqref{H2} follows that the matrix $\hp_k$ is necessarily invertible for all $k\in\di$, and moreover the identities $(\hp_k\pm\hq_k)^{-1}=\hp_k^T\mp\hq_k^T$ and $(\hp_k^{-1}\hq_k)^T=\hp_k^{-1}\hq_k$ hold, see \cite[identities~(11)~and~(12)]{oD.zP01}.
\end{rem}

The following definition was not stated in \cite{oD.zP01} explicitly, but it corresponds to the system from \cite[identity~(19)]{oD.zP01}.
\begin{dfn}[Discrete hyperbolic system]
\label{D4.1}
The system
\begin{equation}
X_\kj=\hp_k X_k+\hq_k U_k, \quad U_\kj=\hq_k X_k+\hp_k U_k, \label{H} \tag{DHS}
\end{equation}
where the matrices $\hp_k$ and $\hq_k$ satisfy identities \eqref{H1} and \eqref{H2} for all $k\in\di$, is called a~\emph{discrete hyperbolic system} \eqref{H}.
\end{dfn}

\begin{rem}
\label{4.3}
Similarly to Remark~\ref{3.3} we can show the justifiability of the terminology ``hyperbolic'' system. It follows from \cite[formulaes~(27)~and~(28)]{oD.zP01}, that the solution of the scalar discrete hyperbolic problem
\begin{equation*}
x_\kj=p_kx_k+q_ku_k, \quad \text{and} \quad u_\kj=q_kx_k+p_ku_k,
\end{equation*}
with the initial conditions $x_0=0$ and $u_0=1$ can be expressed in the form
\[
x_k=\left(\prod\limits_{i=0}^{k-1}\sgn p_i\right)\sinh\!\left(\sum\limits_{i=0}^{k-1}\ln\abs{p_i+q_i}\right)
\quad \text{and} \quad
u_k=\left(\prod\limits_{i=0}^{k-1}\sgn p_i\right)\cosh\!\left(\sum\limits_{i=0}^{k-1}\ln\abs{p_i+q_i}\right).
\]
\end{rem}

By analogy to Lemma~\ref{3.4} we may proove the following lemma.
\begin{lmm}
\label{L4.1}
The pair $\smatrix{X\\U}$ solves the discrete hyperbolic system \eqref{H} if and only if the pair $\smatrix{U\\X}$ solves the same hyperbolic system.
\end{lmm}

The following definition was published in \cite[Definition~3.1]{oD.zP01}.
\begin{dfn}[Discrete matrix-valued hyperbolic functions]
\label{D4.2}
Let $\kn\in\dii$ be fixed. We define the discrete $n \times n$ matrix-valued functions \emph{hyperbolic sine} (denoted by $\sh_{k;\kn}$) and \emph{hyperbolic cosine} (denoted $\ch_{k;\kn}$) by
\[
\sh_{k;\kn}\!:\,=X_k, \quad \ch_{k;\kn}\!:\,=U_k,
\]
where the pair $\smatrix{X\\U}$ is the principal solution of system \eqref{H} at $\kn$. We suppress the index $\kn$ when $\kn=0$, i.e., we denote $\sh_k\!:\,=\sh_{k;0}$ and $\ch_k\!:\,=\ch_{k;0}$.
\end{dfn}

\begin{rem}
\label{P4.12.1}
The solutions of the scalar continuous time system \eqref{CHS} take the form $\sh(t)=\sinh\int_a^t\hq(\tau)\,{\mathrm d}\tau$ and $\ch(t)=\cosh\int_a^t\hq(\tau)\,{\mathrm d}\tau$, see \cite[pg.~12]{kF87}. Hence, we can see discrete matrix functions $\sh_{k;\kn}$ and $\ch_{k;\kn}$ as $n$-dimensional discrete analogs of the scalar continuous hyperbolic functions $\sinh(t-s)$ and $\cosh(t-s)$, which are solutions of the scalar system \eqref{CHS} with $\hq(t)\equiv1$.
\end{rem}

Since the solutions $\smatrix{\ch\\\sh}$ and $\smatrix{\sh\\\ch}$ form normalized conjoined bases of \eqref{H}, it follows that
\begin{equation*}
\label{mphih}
\Psi_k\!:\,=\begin{pmatrix}
\ch_k & \sh_k\\
\sh_k & \ch_k
\end{pmatrix}
\end{equation*}
is a~fundamental matrix of system \eqref{H}. Therefore, every solution $\smatrix{X\\U}$ of \eqref{H} has the form
\[
X_k=\ch_kX_0+\sh_kU_0 \quad \text{and} \quad U_k=\sh_kX_0+\ch_kU_0
\]
for all $k\in\dii$, see also \cite[Theorem~3.1]{oD.zP01}.

From formulas \eqref{5}--\eqref{10} we get the following properties for solutions of discrete hyperbolic systems, which correspond to \cite[Corollary~3.1~and~Theorem~3.2]{oD.zP01}.
\begin{col}
\label{V4.3}
For all $k \in \dii$ the identities
\begin{align}
\ch_k^{T}\ch_k-\sh_k^{T}\sh_k=I&=\ch_k\ch_k^{T}-\sh_k\sh_k^{T}, \label{4.39}\\
\ch^{T}_k\sh_k-\sh^{T}_k\ch_k=0&=\ch_k\sh_k^{T}-\sh_k\ch_k^{T} \label{4.40}
\end{align}
hold, while for all $t \in \di$ we have the identities
\begin{align*}
\ch_\kj\ch_k^{T}-&\sh_\kj\sh^{T}_k=\hp_k,\\
\sh_\kj\ch^{T}_k-&\ch_\kj\sh^{T}_k=\hq_k.
\end{align*}
\end{col}

Similarly to Corollary~\ref{3.8}, with using the Frobenius norm, we can establish a matrix analog of the formula $\cosh^2(t)-\sinh^2(t)=1$, which is \cite[Corollary~3.2]{oD.zP01}.
\begin{col}
\label{D4.4}
For all $k \in \dii$ the next identity holds
\begin{equation*}
\label{4.43}
\norm{\ch}_{F}^{2}-\norm{\sh}_{F}^{2}=n.
\end{equation*}
\end{col}

\begin{rem}
\label{P4.10}
The matrix-valued function hyperbolic cosine has the natural property of scalar hyperbolic cosine. It follows from \eqref{4.39} that the matrix $\ch_k$ is invertible for all $k\in\dii$. Moreover, from \eqref{4.39} and \eqref{4.40} we get
\begin{equation*}
\label{4.72}
\ch_k^{-1}=\ch_k^{T}-\sh_k^{T}\ch_k^{T-1}\sh_k^{T}.
\end{equation*}
On the other hand, the invertibility of the matrix $\sh_k$ is not automatically guaranteed. However, if the matrix $\sh$ is invertible at some point $k\in\dii$, then
\begin{equation*}
\label{4.73}
\sh_k^{-1}=\ch_k^{T}\sh_k^{T-1}\ch_k^{T}-\sh^{T}_k.
\end{equation*}
\end{rem}

The following sum and difference formulas were established in \cite[Theorem~3.3~and~Corollary~3.3]{oD.zP01}.
\begin{thm}
\label{V4.5}
For $k,l\in\dii$ we have
\begin{align}
\sh_{k;l}&=\sh_k\ch^{T}_l-\ch_k\sh^{T}_l, \label{4.44.1}\\
\ch_{k;l}&=\ch_k\ch^{T}_l-\sh_k\sh^{T}_l,\label{4.44.2} \\
\sh_k&=\sh_{k;l}\ch_l+\ch_{k;l}\sh_l,\label{4.45.1}\\
\ch_k&=\ch_{k;l}\ch_l+\sh_{k;l}\sh_l. \label{4.45.2}
\end{align}
\end{thm}

\begin{rem}
\label{P4.5}
With respect to Remark \ref{P4.12.1} for the scalar continuous time case with $\hq(t)\equiv 1$, identities \eqref{4.44.1}--\eqref{4.44.2} are matrix analogues of
\begin{align*}
\sinh(t-s)=\sinh(t)\cosh(s)-\cosh(t)\sinh(s),\\
\cosh(t-s)=\cosh(t)\cosh(s)-\sinh(t)\sinh(s).
\end{align*}
\end{rem}

By interchanging the parameters $k$ and $l$ in \eqref{4.44.1} and \eqref{4.44.2} we get the following corollary, which we can see as $n$-dimensional discrete extensions of the statement about the parity for the scalar hyperbolic functions, see also \cite[Corollary~3.4]{oD.zP01}.
\begin{col}
\label{D4.6}
Let be $k,l\in\dii$. Then the identities
\begin{equation}
\label{4.46}
\sh_{k;l}=-\sh^{T}_{l;k} \quad \text{and} \quad \ch_{k;l}=\ch^{T}_{l;k}
\end{equation}
are true.
\end{col}

Now, by using the same idea and technique as for the discrete trigonometric functions, we can derive analogous results to Theorem~\ref{3.14} for discrete hyperbolic systems. For continuous time hyperbolic systems it was shown in \cite[Theorem~4.2]{kF87}. Hence, we consider the following two discrete hyperbolic systems
\begin{equation}
\label{Hi}
X_\kj=\hp^{(i)}_kX_k+\hq^{(i)}_kU_k,\quad
U_\kj=\hq^{(i)}_kX_k+\hp^{(i)}_kU_k
\end{equation}
with initial conditions $X^{(i)}_0=0$ and $U^{(i)}_0=I$, where $i=1,2$. Denote by $\sh^{(i)}_k$ and $\ch^{(i)}_k$ the corresponding matrix hyperbolic sine and hyperbolic cosine functions from Definition~\ref{D4.2}. Note, that for simplicity we use the following notation $(\hq_k^{(i)})^T=\hq^{(i)T}_k$, $(\hp^{(i)}_k)^T=\hp^{(i)T}_k$, $(\sh^{(i)}_k)^T=\sh^{(i)T}_k$, $(\ch^{(i)}_k)^T=\ch^{(i)T}_k$, $(\shp_k)^T=\sh^{\!+T}_k$, $(\chp_k)^T=\ch^{\!+T}_k$, $(\shm_k)^T=\sh^{\!-T}_k$, and $(\chm_k)^T=\ch^{\!-T}_k$. Now, we put
\begin{align}
\shp_k&\!:\,=\sho_k\ch^{(2)T}_k+\cho_k\sh^{(2)T}_k, \label{4.47}\\
\chp_k&\!:\,=\cho_k\ch^{(2)T}_k+\sho_k\sh^{(2)T}_k, \label{4.48}\\
\shm_k&\!:\,=\sho_k\ch^{(2)T}_k-\cho_k\sh^{(2)T}_k, \label{4.49}\\
\chm_k&\!:\,=\cho_k\ch^{(2)T}_k-\sho_k\sh^{(2)T}_k. \label{4.50}
\end{align}

For two different discrete hyperbolic systems with the same initial conditions the following sum and difference ``formulaes'' hold. For two same discrete hyperbolic systems with different initial conditions it is shown in Theorem~\ref{V4.5}.
\begin{thm}
\label{V4.7}
Assume that $\hp^{(i)}_k$ and $\hq^{(i)}_k$, $i=1,2$, satisfy \eqref{H1} and \eqref{H2}. The pair $\shp_k$ and $\chp_k$ solves the system
\begin{align*}
X_{k+1}&=\hpo_k(X_k\hp^{(2)T}_k+U_k\hq^{(2)T}_k)+\hqo_k(X_k\hq^{(2)T}_k+U_k\hp^{(2)T}_k),\\
U_{k+1}&=\hqo_k(X_k\hp^{(2)T}_k+U_k\hq^{(2)T}_k)+\hpo_k(X_k\hq^{(2)T}_k+U_k\hp^{(2)T}_k),
\end{align*}
with the initial conditions $X_0=0$ and $U_0=I$. The pair $\shm_k$ and $\chm_k$ solves the system
\begin{align*}
X_{k+1}&=\hpo_k(X_k\hp^{(2)T}_k-U_k\hq^{(2)T}_k)+\hqo(-X_k\hq^{(2)T}_k+U_k\hp^{(2)T}_k),\\
U_{k+1}&=\hqo_k(X_k\hp^{(2)T}_k-U_k\hq^{(2)T}_k)+\hpo(-X_k\hq^{(2)T}_k+U_k\hp^{(2)T}_k).
\end{align*}
with the initial conditions $X_0=0$ and $U_0=I$. Moreover, for all $k\in\dii$ we have
\begin{align}
\chp_k\ch^{\!+T}_k-\shp_k\sh^{\!+T}_k=I&=\ch_k^{\!+T}\chp_k-\sh_k^{\!+T}\shp_k, \label{4.51}\\
\chm_k\ch^{\!-T}_k-\shm_k\sh^{\!-T}_k=I&=\ch_k^{\!-T}\chm_k-\sh_k^{\!-T}\shm_k, \label{4.52}\\
\shp_k\ch^{\!+T}_k-\chp_k\sh^{\!+T}_k=0&=\sh_k^{\!+T}\chp_k-\ch_k^{\!+T}\shp_k, \label{4.53}\\
\shm_k\ch^{\!-T}_k-\chm_k\sh^{\!-T}_k=0&=\sh_k^{\!-T}\chm_k-\ch_k^{\!-T}\shm_k. \label{4.54}
\end{align}
\end{thm}

\begin{rem}
\label{P4.4.1}
An analogous statement as in Remark~\ref{3.15} for the pairs $\smatrix{\stp\\ \ctp}$ and $\smatrix{\stm\\ \ctm}$ now holds for the solutions $\smatrix{\shp\\ \chp}$ and $\smatrix{\shm\\ \chm}$. It means, these two pairs are not conjoined bases of their corresponding systems, because these systems are not symplectic.
\end{rem}

As in Corollary~\ref{3.16}, when the two systems in \eqref{Hi} are the same, we obtain from Theorem~\ref{V4.7} the following. It is a discrete analogue of \cite[Corollary~1]{kF87}.
\begin{col}
\label{D4.8}
Assume that $\hp_k$ and $\hq_k$ satisfy \eqref{H1} and \eqref{H2}. Then the system
\begin{align*}
X_{k+1}&=\hp_k(X_k\hp_k^T+U_k\hq_k^T)+\hq_k(X_k\hq_k^T+U_k\hp_k^T),\\
U_{k+1}&=\hq_k(X_k\hp_k^T+U_k\hq_k^T)+\hp_k(X_k\hq_k^T+U_k\hp_k^T)
\end{align*}
with the initial conditions $X_0=0$ and $U_0=I$ possesses the solution
\[
X_k=2\sh_k\ch_k^{T} \quad \text{and} \quad U=\ch_k\ch_k^{T}+\sh_k\sh_k^{T},
\]
where $\sh$ and $\ch$ are the matrix functions in Definition~\ref{D4.2}. Moreover, the above matrices $X$ and $U$ commute, i.e. $X_kU_k=U_kX_k$.
\end{col}

\begin{rem}
\label{P4.6.1}
The previous corollary can be viewed as the $n-$dimensional discrete analogy of the double angle formulae for scalar continuous time hyperbolic functions
\[
\sinh(2t)=2\sinh(t)\cosh(t) \quad \text{and} \quad \cosh(2t)=\cosh^{2}(t)+\sinh^{2}(t).
\]
\end{rem}

Now we can formulate similar identities as for trigonometric functions in Corollary~\ref{3.18}.
\begin{col}
\label{D4.9}
For all $k \in \dii$ we have the identities
\begin{align}
\sho_k\sh_k^{(2)T}&=\frac{1}{2}\,(\chp_k-\chm_k),\label{4.56}\\
\cho_k\ch_k^{(2)T}&=\frac{1}{2}\,(\chp_k+\chm_k),\label{4.57}\\
\sho_k\ch_k^{(2)T}&=\frac{1}{2}\,(\shp_k+\shm_k).\label{4.58}
\end{align}
\end{col}

\begin{rem}
\label{P4.9}
Identities \eqref{4.56}--\eqref{4.58} are discrete $n$-dimensional analogues of
\begin{align*}
\sinh(t)\sinh(s)&=\frac{1}{2}\,[\cosh(t+s)-\cosh(t-s)],\\
\cosh(t)\cosh(s)&=\frac{1}{2}\,[\cosh(t+s)+\cosh(t-s)],\\
\sinh(t)\cosh(s)&=\frac{1}{2}\,[\sinh(t+s)+\sinh(t-s)].
\end{align*}
\end{rem}

The next definition of discrete matrix hyperbolic functions tangent and cotangent is known from \cite[Definition~3.2]{oD.zP01}. Recall that the matrix function $\ch_k$ is invertible for all $k\in\dii$, see Remark~\ref{P4.10}.
\begin{dfn}
\label{D4.3}
We define the discrete matrix-valued function \emph{hyperbolic tangent} (we write $\htg$) by
\[
\htg_k\!:\,=\ch^{-1}_k\,\sh_k.
\]
Whenever $\sh_k$ is invertible we define the discrete matrix-valued function \emph{hyperbolic cotangent} (we write $\hco$) by
\[
\hco_k\!:\,=\sh^{-1}_k\,\ch_k.
\]
\end{dfn}

Analogous results for discrete hyperbolic tangent and cotangent as in Theorem~\ref{3.21} and Theorem~\ref{3.22} for trigonometric functions tangent and cotangent were published in \cite[Theorem~3.4]{oD.zP01}.
\begin{thm}
\label{V4.11}
For $k\in\dii$ we have
\begin{align}
\htg^{T}_k&=\htg_k,\label{4.68}\\
\ch^{-1}_k\ch^{T-1}&_k+\htg^{2}_k=I, \label{4.69}
\end{align}
and for $k\in\di$
\begin{equation}
\label{4.70}
\Delta\htg_k=\ch_\kj^{-1}\hq_k\,\ch^{T-1}_k.
\end{equation}
\end{thm}

\begin{rem}
\label{P4.11}
In the scalar case $n=1$ identity \eqref{4.68} is trivial. In the scalar continuous time case \eqref{CHS} identity \eqref{4.69} takes the form
\[
\frac{1}{\cosh^{2}(s)}+\tanh^{2}(s)=1 \quad \text{with \ $s=\int_a^t\!\!\!\hq(\tau)\,{\mathrm d}\tau$.}
\]
Finally, identity \eqref{4.70} takes the form
\[
\left(\tanh\!\int_a^t\!\!\!\hq(\tau)\,{\mathrm d}\tau\right)'=\frac{\hq(t)}{\cosh^{2}\int_a^t\hq(\tau)\,{\mathrm d}\tau}.
\]
\end{rem}

\begin{thm}
\label{V4.12}
Whenever $\hco_k$ is defined we get
\begin{align}
\hco^{T}_k&=\hco_k,\label{4.76}\\
\hco^{2}_k-\sh&^{-1}_k\,\sh^{T-1}_k=I. \label{4.77}
\end{align}
Moreover, if\, $\sh_k$ and $\sh_\kj$ are invertible, then
\begin{equation}
\label{4.78}
\Delta\hco_k=-\sh_\kj^{-1}\hq_k\,\sh^{T-1}_k.
\end{equation}
\end{thm}

\begin{rem}
\label{P4.12}
In the scalar case $n=1$ identity \eqref{4.76} is trivial. In the scalar continuous time case \eqref{CHS} identity \eqref{4.77} reads as
\[
\cotanh^{2}(s)-\frac{1}{\sinh^{2}(s)}=1,
\]
where $s\neq0$ is from Remark \ref{P4.11}. Finally, for $\int_a^t\hq(\tau)\,{\mathrm d}\tau\neq0$ reduces identity \eqref{4.78} to
\[
\left(\cotanh\!\int_a^t\!\!\!\hq(\tau)\,{\mathrm d}\tau\right)'=\frac{-\hq(t)}{\sinh^{2}\int_a^t\hq(\tau)\,{\mathrm d}\tau}.
\]
\end{rem}

Next, similarly to the definitions of the time scale matrix functions $\sh^{(i)}_k$, $\ch^{(i)}_k$, $\sh^{\pm}_k$, and $\ch^{\pm}_k$ from \eqref{Hi}--\eqref{4.50} we define
\begin{align*}
\htg^{(i)}_k&\!:\,=(\ch^{(i)}_k)^{-1}\,\sh^{(i)}_k,&\quad\hco^{(i)}_k&\!:\,=(\sh^{(i)}_k)^{-1}\,\ch^{(i)}_k,\\
\htgp_k&\!:\,=(\chp_k)^{-1}\shp_k,&\quad\hcop_k&\!:\,=(\shp_k)^{-1}\chp_k,\\
\htgm_k&\!:\,=(\chm_k)^{-1}\shm_k,&\quad\hcom_k&\!:\,=(\shm_k)^{-1}\chm_k.
\end{align*}

\begin{rem}
\label{P4.14}
As in Remark~\ref{3.23} we conclude that the first identities from \eqref{4.53} and \eqref{4.54} imply the symmetry of the functions $\htg^{\pm}_k$ and $\hco^{\pm}_k$.
\end{rem}

\begin{thm}
\label{V4.13}
For all $k\in\dii$  such that all involved functions are defined we have
\begin{align}
\htgo_k+\htgt_k&=\htgo_k\,(\hcoo_k+\hcot_k)\,\htgt_k, \label{4.84.t}\\
\htgo_k-\htgt_k&=\htgo_k\,(\hcot_k-\hcoo_k)\,\htgt_k, \label{4.85.t}\\
\htgo_k+\htgt_k&=(\cht_k)^{-1}\,\sh_k^{\!+T}(\cho_k)^{T-1}, \label{4.79}\\
\htgo_k-\htgt_k&=(\cht_k)^{-1}\,\sh_k^{\!-T}(\cho_k)^{T-1}, \label{4.80}\\
\htgp_k&=(\cht_k)^{T-1}\,(I+\htgo_k\htgt_k)^{-1}\, \notag \\
  &\ \hspace*{45mm} \times(\htgo_k+\htgt_k)\,\ch^{(2)T}_k, \label{4.82}\\
\htgm_k&=(\cht_k)^{T-1}\,(I-\htgo_k\htgt_k)^{-1}\, \notag \\
  &\ \hspace*{45mm} \times(\htgo_k-\htgt_k)\,\ch^{(2)T}_k, \label{4.83}\\
\hcoo_k+\hcot_k&=\hcoo_k\,(\htgo_k+\htgt_k)\,\hcot_k, \label{4.84}\\
\hcoo_k-\hcot_k&=\hcoo_k(\htgt_k-\htgo_k)\,\hcot_k, \label{4.85}\\
\hcoo_k+\hcot_k&=(\sht_k)^{-1}\,\sh_k^{\!+T}(\sho_k)^{T-1}, \label{4.86}\\
\hcoo_k-\hcot_k&=-(\sht_k)^{-1}\,\sh_k^{\!-T}(\sho_k)^{T-1}. \label{4.87}\\
\hcop_k&=(\sht_k)^{T-1}\,(\hcoo_k+\hcot_k)^{-1}\, \notag \\
  &\ \hspace*{40mm} \times(\hcoo_k\hcot_k+I)\,\sh^{(2)T}_k, \label{4.82.c}\\
\hcom_k&=(\sht_k)^{T-1}\,(\hcot_k-\hcoo_k)^{-1}\, \notag \\
  &\ \hspace*{40mm} \times(\hcoo_k\hcot_k-I)\,\sh^{(2)T}_k, \label{4.83.c}
\end{align}
\end{thm}

\begin{rem}
\label{P4.15}
Consider now the scalar continuous time case \eqref{CHS} with $\hq(t)\equiv 1$. Then identities  \eqref{4.84.t} and \eqref{4.85.t} are equivalent to
\begin{align*}
\tanh(t)+\tanh(s)&=[\cotanh(t)+\cotanh(s)]\tanh(t)\tanh(s),\\
\tanh(t)-\tanh(s)&=[\cotanh(s)-\cotanh(t)]\tanh(t)\tanh(s),
\end{align*}
respectively, identities \eqref{4.79} and \eqref{4.80} have the form
\[
\tanh(t)+\tanh(s)=\frac{\sinh(t+s)}{\cosh(t)\cosh(s)},\quad
\tanh(t)-\tanh(s)=\frac{\sinh(t-s)}{\cosh(t)\cosh(s)},
\]
respectively, while identities \eqref{4.82} and \eqref{4.83} have the form
\[
\tanh(t+s)=\frac{\tanh(t)+\tanh(s)}{1+\tanh(t)\tanh(s)},\quad
\tanh(t-s)=\frac{\tanh(t)-\tanh(s)}{1-\tanh(t)\tanh(s)}.
\]
respectively. Moreover, the identities in \eqref{4.84} and \eqref{4.85} reduce to
\begin{align*}
\cotanh(t)+\cotanh(s)=[\tanh(t)+\tanh(s)]\cotanh(t)\cotanh(s),\\
\cotanh(t)-\cotanh(s)=[\tanh(s)-\tanh(t)]\cotanh(t)\cotanh(s),
\end{align*}
respectively, and the identities in \eqref{4.86} and \eqref{4.87} correspond in this case to
\[
\cotanh(t)+\cotanh(s)=\frac{\sinh(t+s)}{\sinh(t)\,\sinh(s)},\quad
\cotanh(t)-\cotanh(s)=\frac{\sinh(s-t)}{\sinh(t)\,\sinh(s)},
\]
respectively. Finally, the identities \eqref{4.82.c} and \eqref{4.83.c} is in accordance with
\[
\cotanh(t+s)=\frac{\cotanh(t)\cotanh(s)+1}{\cotanh(t)+\cotanh(s)},\quad
\cotanh(t-s)=\frac{\cotanh(t)\cotanh(s)-1}{\cotanh(s)-\cotanh(t)},
\]
respectively.
\end{rem}
\section{Concluding remarks}
\label{sec:5}

In this paper we extended to discrete matrix case some identities known for the scalar continuous time trigonometric and hyperbolic functions (for an overview see e.g. \cite[Chapter~4]{mA.iaS64}). On the other hand, there are still several trigonometric and hyperbolic identities which we could not extend to the discrete $n$-dimensional case, e.g.
\begin{align}
\left.\begin{alignedat}{2}
\label{5.1}
\sin{x}\pm\sin{y}&=2\sin{\frac{x\pm y}{2}}\cos{\frac{x\mp y}{2}},\\
\cos{x}+\cos{y}&=2\cos{\frac{x+y}{2}}\cos{\frac{x-y}{2}},\\
\cos{x}-\cos{y}&=2\sin{\frac{x+y}{2}}\sin{\frac{y-x}{2}},\\
\sinh{x}\pm\sinh{y}&=2\sinh{\frac{x\pm y}{2}}\cosh{\frac{x\mp y}{2}},\\
\cosh{x}+\cosh{y}&=2\cosh{\frac{x+y}{2}}\cosh{\frac{x-y}{2}},\\
\cosh{x}-\cosh{y}&=2\sinh{\frac{x+y}{2}}\sinh{\frac{x-y}{2}},\\
\end{alignedat}\right\}
\end{align}
\begin{align}
\left.\begin{alignedat}{2}
\label{5.2}
\sin{(x+y)}\sin{(x-y)}&=\sin^{2}{x}-\sin^{2}{y},\\
\cos{(x+y)}\cos{(x-y)}&=\cos^{2}{x}-\sin^{2}{y},\\
\sinh{(x+y)}\sinh{(x-y)}&=\sinh^{2}{x}-\sinh^{2}{y},\\
\cosh{(x+y)}\cosh{(x-y)}&=\sinh^{2}{x}+\cosh^{2}{y}.
\end{alignedat}\right\}
\end{align}

The first identity in \eqref{5.1} reduces to $\sin{x}=2\sin{\frac{x}{2}}\cos{\frac{x}{2}}$ if we put $y=0$. Its right-hand side seems to be similar to the solution $X_k$ of the system in Corollary~\ref{3.16} but the left-hand side is not the matrix-valued function $\st_k$ because this system is not trigonometric. Similar reasoning can be used for the remaining identities in \eqref{5.1}.

Nevertheless, we can calculate the corresponding products for matrix-valued functions $\stp_k$, $\ctp_k$, $\shp_k$, and $\chp_k$ as on left-hand side in \eqref{5.2}. But we can not to modify these products in to the forms, which are similar to the right-hand side in \eqref{5.2}, because the matrix multiplication is not commutative.

\section{Acknowledgement}
\label{sec:6}

The author thanks to his advisor Roman \v{S}imon Hilscher for his help and remarks which led to the present version of the paper.


\end{document}